\documentclass[12pt]{amsart}

\theoremstyle{definition}
\newtheorem{problem}{Problem}

\newcommand{\C}{\mathbb{C}}
\DeclareMathOperator{\re}{Re}
\renewcommand{\iint}{\int\!\!\!\int}   
\newcommand{\Beta}{B}

\begin{document}

\title{Lu Qi-Keng's Problem}

\author{Harold P. Boas}

\address{Department of Mathematics, Texas A\&M University,
  College Station, TX 77843--3368, USA}

\email{boas@math.tamu.edu}

\thanks{This article is based on a lecture at the third Korean
  several complex variables symposium in December 1998.  The
  author thanks GARC, the Global Analysis Research Center at
  Seoul National University, for sponsoring his participation in
  this international conference.}  

\thanks{The author's research was partially supported by grant
  number DMS-9801539 from the National Science Foundation of the
  United States of America.}

\begin{abstract}
  This expository article, intended to be accessible to students,
  surveys results about the presence or absence of zeroes of the
  Bergman kernel function of a bounded domain in~\(\C^n\). Six
  open problems are stated.
\end{abstract}

\subjclass{Primary 32H10}

\maketitle

\setcounter{tocdepth}{1}
\tableofcontents

\section{Introduction}
When does a convergent infinite series of holomorphic functions
have zeroes?  This question is a fundamental, difficult problem
in mathematics.

For example, the series \(\sum_{n=0}^\infty z^n/n!\) is
zero-free, but how can one tell this without a priori knowledge
that the series represents the exponential function? Changing the
initial term of this series produces a new series that does have
zeroes, since the range of the exponential function is all
non-zero complex numbers. The problem of determining when a
series has zeroes is essentially equivalent to the hard problem
of determining the range of a holomorphic function that is
presented as a series.

A famous instance of the problem of locating zeroes of infinite
series is the Riemann hypothesis about the zeta function: namely,
the conjecture that when \(0<\re z<1/2\), the convergent series
\(\sum_{n=1}^\infty (-1)^n/n^z\) has no zeroes. This formulation
of the Riemann hypothesis is equivalent to the usual statement
that the zeroes of \(\zeta(z)\) in the critical strip where
\(0<\re z<1\) all lie on the line where \(\re z=1/2\). Indeed,
when \(\re z>1\), absolute convergence justifies writing that
\begin{equation*}
\sum_{n=1}^\infty \frac{(-1)^n}{n^z} = -\sum_{n=1}^\infty
\frac{1}{n^z} +2\sum_{\text{even \(n\)}}\frac{1}{n^z}
=\zeta(z)(2^{1-z}-1).
\end{equation*}
By the principle of persistence of functional relationships, the
expressions on the outer ends of this equality still agree when
\(0<\re z<1\), so \(\zeta(z)\) and \(\sum_{n=1}^\infty
(-1)^n/n^z\) have the same zeroes in the interior of the critical
strip. Moreover, the well known functional equation for the zeta
function implies that the zeroes of~\(\zeta\) in the critical
strip are symmetric about the point~\(1/2\), so it suffices to
examine the left half of the critical strip for zeroes.

In this article, I shall discuss a different instance of the
general problem of locating zeroes of infinite series.  The
Bergman kernel function is most conveniently expressed as the sum
of a convergent infinite series.  In 1966, Lu Qi-Keng \cite{lu}
asked: for which domains is the Bergman kernel function
\(K(z,w)\) zero-free? I shall address three general methods for
approaching this problem, and I shall give examples both of
domains whose Bergman kernel functions are zero-free and of
domains whose Bergman kernel functions have zeroes.

\section{The Bergman kernel function}
\subsection{Definition}
The Bergman kernel function \(K(z,w)\) of a domain~\(D\) in
\(\C^n\) is the unique
sesqui-holomorphic\footnote{Sesqui-holomorphic means holomorphic
  in the first variable and conjugate holomorphic in the second
  variable.} function satisfying the skew-symmetry property that
\(K(z,w)=\overline {K(w,z)}\) and the reproducing property that
\begin{equation*}
f(z)=\int_D K(z,w) f(w) \,d\text{Volume}_w  \qquad \text{for all
  $z$ in $D$}
\end{equation*}
for every square-integrable holomorphic function~\(f\) on~\(D\).
Equivalently, \(K(z,w)=\sum_j
\varphi_j(z)\overline{\varphi_j(w)}\), where \(\{\varphi_j\}\) is
an orthonormal basis for the Hilbert space of square-integrable
holomorphic functions on~\(D\). To compute the Bergman kernel
function, one typically chooses an orthogonal basis, calculates
the normalizing factors, and sums the series.

For example, the monomials \(1\), \(z\), \(z^2\), 
\ldots, are orthogonal on the unit disk in the complex plane, and
the norm of \(z^k\) is \((\int_{0}^{2\pi} \int_{0}^{1}
r^{2k+1}\,dr\,d\theta)^{1/2}\), or \(\sqrt{2\pi/(2k+2)}\).
Therefore the Bergman kernel function \(K(z,w)\) of the unit disk
equals
\begin{equation}
  \label{eq:disk}
  \sum_{k=0}^\infty \frac{k+1}{\pi} \cdot z^k \bar w^k 
  = \frac{1}{\pi}\cdot \frac{1}{(1-z\bar w)^2}.
\end{equation}
This result is compatible with the Cauchy integral formula
\begin{equation*}
f(z)=\frac{1}{2\pi i}\oint_{|w|=1}
\frac{f(w)}{w-z}\,dw,
\end{equation*}
which can be rewritten by Green's formula as
\begin{equation*}
f(z)=\frac{1}{\pi}\iint_{|w|<1} \frac{f(w)}{(1-z\bar w)^2}
\,d\text{Area}_w,
\end{equation*}
thus confirming that the kernel function~\eqref{eq:disk} does
have the required reproducing property.
 
\subsection{Transformation rule}
If \(f:D_1\to D_2\) is a biholomorphic\footnote{Biholomorphic
  means holomorphic with a holomorphic inverse.} mapping, and if
\(K_1\) and~\(K_2\) denote the Bergman kernel functions of the
domains \(D_1\) and~\(D_2\) in~\(\C^n\), then
\begin{equation}
\label{eq:transform}
K_1(z,w)=(\det f'(z)) K_2(f(z), f(w)) (\,\overline{\det f'(w)}\,).
\end{equation}
This relationship holds because if \(\{\varphi_j\}\) is an
orthonormal basis for the square-integrable holomorphic functions
on~\(D_2\), then \(\{(\det f') \,\varphi_j\circ f\}\) is an
orthonormal basis for the square-integrable holomorphic functions
on~\(D_1\). For example, scaling~\eqref{eq:disk} shows that the
Bergman kernel function of the disk of radius~\(r\) is equal to
\(\dfrac{1}{\pi r^2}\cdot \dfrac{1}{(1-z\bar w/r^2)^2}\).

It is an observation of Steven~R. Bell \cite{bell1981, bell1982}
that a similar transformation rule holds even when \(f:D_1\to
D_2\) is a branched \(m\)-fold covering (a proper holomorphic
mapping): namely,
\begin{equation*}
\sum_{k=1}^m K_1(z, f_k^{-1}(w)) 
(\,\overline{\det (f_k^{-1})'(w)}\,) =
(\det f'(z)) K_2(f(z),w),
\end{equation*}
where the \(f_k^{-1}\) are the \(m\) holomorphic local inverses
of~\(f\).  This formula is not valid on the branching set, where
local inverses are not defined.

The Bergman kernel function \(K(z,w)\) of a simply-connected
planar domain~\(D\) is related to the Riemann mapping
function~\(f\) that maps~\(D\) onto the unit disk, taking the
point~\(a\) to~\(0\): namely,
\begin{equation}
\label{eq:riemann}
f'(z)=K(z,a)\sqrt{\frac{\pi}{K(a,a)}}.
\end{equation}
Indeed, the  transformation rule~\eqref{eq:transform} implies that
\begin{equation*}
K(z,w)=f'(z)\cdot\frac{1}{\pi}\cdot
\frac{1}{(1-f(z)\overline{f(w)}\,)^2}\cdot\overline{f'(w)}.
\end{equation*}
Since \(f(a)=0\), and \(f'(a)\) is real and positive,
setting \(w=a\) implies that
\begin{math}
\pi K(z,a)=f'(z) f'(a)
\end{math},
and then setting \(z=a\) makes it possible to eliminate~\(f'(a)\)
to obtain~\eqref{eq:riemann}.

Since the Riemann mapping function solves a certain extremal
problem, the connection in one dimension with the Bergman kernel
function suggests studying the Bergman kernel function in higher
dimensions in connection with extremal problems.

\subsection{Extremal properties}
\label{extremal}
One can use extremal characterizations of the Bergman kernel
function to help prove theorems about convergence of the Bergman
kernel functions of a convergent sequence of domains.  In the
following two properties, the point~\(w\) is fixed in a
domain~\(D\) in~\(\C^n\), and \(\{\varphi_j\}\) is an orthonormal
basis for the Hilbert space \(A^2(D)\) of square-integrable
holomorphic functions on~\(D\).

\begin{enumerate}
\item In the class of holomorphic functions~\(f\) on~\(D\) such
  that \(\int_D |f|^2\le 1\), the maximal value of \(|f(w)|^2\)
  is \(K(w,w)\). In other words, \(K(w,w)\) is the square of the
  norm of the functional from \(A^2(D)\) to~\(\C\) that
  evaluates a function at the point~\(w\).
  
  Indeed, if \(f(z)=\sum_j c_j \varphi_j(z)\), and if \(\sum_j
  |c_j|^2\le1\), then the Cauchy-Schwarz inequality implies that
  \(|\sum_j c_j \varphi_j(w)|^2\) is bounded above by \(\sum_j
  |\varphi_j(w)|^2\), which equals \(K(w,w)\); and the upper
  bound is attained if \(c_j\) is taken equal to
  \(\overline{\varphi_j(w)}/(\sum_j |\varphi_j(w)|^2)^{1/2}\).
  
\item In the class of holomorphic functions~\(f\) on~\(D\)
  satisfying the nonlinear constraint that \(f(w)\ge \int_D
  |f|^2\), the function with the maximal value at~\(w\) is
  \(K(\cdot,w)\).
  
  Indeed, it is evident that the function \(K(\cdot,w)= \sum_j
  \overline{\varphi_j(w)}\,\varphi_j\) is in the class. On the
  other hand, if \(f=\sum_j c_j \varphi_j\) is an arbitrary
  member of the class, then the preceding extremal property
  implies that \(f(w)^2\le K(w,w)\sum_j |c_j|^2\); and since the
  defining property of the class implies that \((\sum_j
  |c_j|^2)^2 \le f(w)^2\), it follows that \(\sum_j|c_j|^2\le
  K(w,w)\), and hence \(f(w)\le K(w,w)\).
\end{enumerate}

\section{Motivation for Lu Qi-Keng's problem}
The Riemann mapping theorem characterizes the planar domains that
are biholomorphically equivalent to the unit disk. In higher
dimensions, there is no Riemann mapping theorem,\footnote{More
  precisely, in order to obtain a generalized Riemann mapping
  theorem, one needs either new hypotheses \cite{chernji} or new
  definitions \cite{semmes}.} and two natural problems arise.

\begin{enumerate}
\item Are there canonical representatives of biholomorphic
  equivalence classes of domains?
  
\item How can one tell that two particular domains are
biholomorphically inequivalent?
\end{enumerate}

As an approach to the first question, Stefan Bergman introduced
the notion of a ``representative domain'' to which a given domain
may be mapped by ``representative coordinates''.  If \(g_{jk}\)
denotes the Bergman metric
\begin{math}
\dfrac{\partial^2}{\partial z_j\partial\bar  z_k}\log K(z,z)
\end{math},
where \(K\)~is the Bergman kernel function, then the local
representative coordinates based at the point~\(a\) are
\begin{equation*}
\sum_{k=1}^n
g_{kj}^{-1}(a) \left. \frac{\partial}{\partial \bar w_k} \log
\frac{K(z,w)}{K(w,w)} \right|_{w=a}, \qquad j=1, \dotsc, n.
\end{equation*}
These coordinates take \(a\) to~\(0\) and have complex Jacobian
matrix at~\(a\) equal to the identity.

Zeroes of the Bergman kernel function \(K(z,w)\) evidently pose
an obstruction to the global definition of Bergman representative
coordinates. This observation was Lu Qi-Keng's motivation for
asking which domains have zero-free Bergman kernel functions.

On the other hand, if the Bergman kernel function of a domain
does have zeroes, then the transformation
rule~\eqref{eq:transform} shows that the zero set is a
biholomorphically invariant object. Therefore zero sets of
Bergman kernel functions could be a tool for addressing the
second question stated above. This idea has not yet been
exploited in the literature.

\section{First examples}
\label{first}
The Bergman kernel function~\eqref{eq:disk} of the unit disk is
evidently zero-free.  Consequently, the Bergman kernel function
of every bounded, simply connected, planar domain is zero-free:
apply either the transformation rule~\eqref{eq:transform} or the
explicit formula~\eqref{eq:riemann} relating the Riemann mapping
function to the Bergman kernel function.

On the other hand, the Bergman kernel function of every annulus
does have zeroes \cite{rosenthal,skwar}; more generally, the
Bergman kernel function of every bounded, multiply connected,
planar domain with smooth boundary has zeroes \cite{suitayamada}.

Isolated singularities of square-integrable holomorphic functions
are removable, and therefore the Bergman kernel function does not
see isolated punctures in a domain.  For example, the Bergman
kernel function of a punctured disk is zero-free. On the other
hand, a finitely connected planar domain with no singleton
boundary component can be mapped biholomorphically to a smoothly
bounded domain.  Consequently, if a bounded planar domain is
finitely connected and has at least two non-singleton boundary
components, then its Bergman kernel function has zeroes.

I do not know if a corresponding statement holds for infinitely
connected planar domains. For example, delete from the open unit
disk a countable sequence of pairwise disjoint closed disks that
accumulate only at the boundary of the unit disk. Does the
Bergman kernel function of the resulting domain have zeroes?

\begin{problem}
  Give necessary and sufficient conditions on an infinitely
  connected planar domain for its Bergman kernel function to have
  zeroes.
\end{problem}

It is easy to see that in higher dimensions, the Bergman kernel
function of a product domain is the product of the Bergman kernel
functions of the lower dimensional domains. Consequently, the
Bergman kernel function of a polydisc is zero-free, while the
Bergman kernel function of the Cartesian product of a disc with
an annulus does have zeroes. The Bergman kernel function
\(K(z,w)\) of the unit ball in~\(\C^n\) is the zero-free function
\(\dfrac{n!}{\pi^n}\cdot \dfrac{1}{(1-\langle
z,w\rangle)^{n+1}}\), where \(\langle z,w\rangle\) denotes the
scalar product \(z_1\bar w_1 + \dots + z_n \bar w_n\).

Even without knowing this explicit formula, one can see that the
Bergman kernel function of the unit ball is zero-free.  Since the
ball is a complete circular domain,\footnote{A domain is called
  complete circular if whenever it contains a point~\(z\), it
  also contains the one-dimensional disk \(\{\,\lambda z:
  |\lambda|\le1\,\}\).}  its space of square-integrable
holomorphic functions has an orthonormal basis whose first
element is a constant (namely, the reciprocal of the square root
of the volume of the domain) and whose other elements are
functions that vanish at the origin.  Consequently, the Bergman
kernel function~\(K\) has the property that \(K(z,0)\) is a
non-zero constant function of~\(z\). Since the ball is
homogeneous,\footnote{A domain is called homogeneous if it has a
  transitive automorphism group, that is, if any point of the
  domain can be mapped to any other point by a biholomorphic
  self-mapping of the domain.}  the transformation
rule~\eqref{eq:transform} implies that the Bergman kernel
function is nowhere zero. The same argument shows that every
bounded, homogeneous, complete circular domain has a zero-free
Bergman kernel function \cite{belldavid}.

For many years it was thought that sufficiently nice,
topologically trivial, bounded domains in~\(\C^n\) should have
zero-free Bergman kernel functions. For example, all strongly
convex,\footnote{The statement is also true for strongly
  pseudoconvex domains, that is, domains that locally can be
  mapped biholomorphically to strongly convex domains.}
sufficiently small perturbations of the ball have zero-free
Bergman kernel functions if ``small'' is interpreted in the
\(C^\infty\) topology on domains \cite{greenekrantz}.  It turns
out, however, to be the generic situation for the Bergman kernel
function of a domain of holomorphy to have zeroes, if ``generic''
is interpreted in the very flexible Hausdorff topology on domains
\cite{boasgeneric}.

\begin{problem}
  Do there exist arbitrarily small class \(C^1\) perturbations of
  the ball whose Bergman kernel functions have zeroes?
\end{problem}

I shall now discuss three techniques that can be used to show the
existence of interesting domains whose Bergman kernel functions
have zeroes.

\section{Variation of domains}
\label{variation-domains}
If reasonable domains \(\Omega_j\) converge in a reasonable way
to a limiting domain \(\Omega\), then the Bergman kernel
functions \(K_{\Omega_j}(z,w)\) converge to \(K_\Omega(z,w)\)
uniformly on compact subsets of \(\Omega\times\Omega\).  The word
``reasonable'' can be made precise \cite{boasgeneric}, but here I
will simply mention two examples of reasonable behavior. The
first example is a fundamental theorem of I.~P. Ramadanov
\cite{ramadanov} which started the whole theory.
\begin{itemize}
\item The \(\Omega_j\) form an increasing sequence whose union
  is~\(\Omega\).
  
\item The \(\Omega_j\) are bounded pseudoconvex\footnote{A domain
is pseudoconvex if it is the union of an increasing sequence of
strongly pseudoconvex domains, as defined in the preceding
footnote. According to the solution of the Levi problem (see, for
example, \cite{Krantz}), pseudoconvex domains are the same as
domains of holomorphy.} domains whose complements converge in the
Hausdorff metric to the complement of an~\(\Omega\) whose
boundary is locally a graph.
\end{itemize}
An example of unreasonable convergence is a sequence of disks
shrinking down to a disk with a slit \cite{skwarczynski}.

The proof of the convergence theorem exploits the extremal
characterization of the Bergman kernel function from
section~\ref{extremal}. The application to Lu Qi-Keng's problem
is that by Hurwitz's theorem, if the Bergman kernel function of
the limiting domain~\(\Omega\) has zeroes, then so does the
Bergman kernel function of the approximating domain~\(\Omega_j\)
when \(j\)~is sufficiently large.

Consequently, to construct a nice domain whose Bergman kernel
function has zeroes, it suffices to construct a degenerate domain
whose Bergman kernel function has zeroes, and then to approximate
the degenerate domain by nice ones. For example, an easy
calculation shows that constant functions are not
square-integrable on the domain
\begin{equation}
\label{eq:counterex}
  \left\{\, (z_1,z_2)\in\mathbf{C}^2: |z_2|<\frac{1}{1+|z_1|}\,
\right\},
\end{equation}
although the domain does support some non-constant
square-integrable holomorphic functions. Therefore the Bergman
kernel function of this domain has a zero at the origin, but is
not identically zero.  By approximating this domain from inside,
one sees that there exists a bounded, smooth, logarithmically
convex, complete Reinhardt domain\footnote{A domain is called
  complete Reinhardt if whenever it contains a point
  \((z_1,\dotsc,z_n)\), it also contains the polydisc
  \(\{\,(\lambda_1z_1, \dotsc, \lambda_nz_n): |\lambda_1|\le1,
  \dotsc, |\lambda_n|\le1\,\}\). The logarithmically convex
  complete Reinhardt domains are the convergence domains of power
  series.}  whose Bergman kernel function has zeroes
\cite{boascounterexample}.  This example was surprising when it
was first discovered; indeed, it contradicts a theorem previously
published by two different authors \cite[Theorem~1]{matsuura},
\cite[Corollary]{kakurai} and applied by a third \cite{kanemaru}.

Nguy\^en Vi\^et Anh, a student in Marseille, recently showed
\cite{anh} how to approximate the domain~\eqref{eq:counterex}
from inside by \emph{concrete} domains that are smooth,
algebraic, logarithmically convex, complete Reinhardt domains.
Namely, the domain defined by the inequality
\begin{equation}
\label{eq:anh}
  |z_2|^{2k}(1+|z_1|)^{2k} + |z_2|^{2k}(1-|z_1|)^{2k} + \left(
  \frac{|z_1|^2+|z_2|^2} {k}\right)^k <1
\end{equation}
has the indicated properties when \(k\)~is a positive integer,
and so the Bergman kernel function of this domain must have
zeroes when \(k\)~is sufficiently large.

\begin{problem}
  How large must \(k\) be in order for the Bergman kernel
  function of the domain~\eqref{eq:anh} to have zeroes?
\end{problem}

The technique of variation of domains can be used to prove the
statement at the end of section~\ref{first} that every nice
domain can be arbitrarily closely approximated in the Hausdorff
metric by a nice domain whose Bergman kernel function has zeroes.
View the starting domain as the Earth, and place in orbit around
the Earth a small copy of one of the bounded domains just
discussed. The Bergman kernel function of this disconnected
region has zeroes because the Bergman kernel function of the
satellite has zeroes. Now attach the satellite to the Earth by a
thin tether. If the tether is sufficiently thin, then the
principle of variation of domains implies that the Bergman kernel
function of the joined domain has zeroes. According to the
barbell lemma \cite[Chapter~5, Exercise~21]{Krantz} in a suitable
formulation \cite{boasgeneric}, the joined domain can be made
strongly pseudoconvex if the Earth is.

\section{Variation of weights}
In the representation of the Bergman kernel function as an
infinite series \(\sum_j \varphi_j(z)\overline{\varphi_j(w)}\),
the basis elements~\(\varphi_j\) are supposed to be orthonormal
for integration with respect to Lebesgue measure.  It is natural
to consider the analogous construction when Lebesgue measure is
multiplied by a positive weight function.

For instance, if \(G\)~is a domain in~\(\C^n\), and \(\Omega\)~is
a Hartogs domain in~\(\C^{n+1}\) with base~\(G\), which means
that \(\Omega = \{\,(z,z_{n+1})\in\C^{n+1}: z\in G\) and
\(|z_{n+1}|<r(z)\,\}\), where \(r\)~is a positive function
on~\(G\), then the Bergman kernel function of~\(\Omega\)
restricted to the base~\(G\) equals the weighted Bergman kernel
function of~\(G\) corresponding to the weight \(\pi r^2\). This
follows because the Bergman kernel function is uniquely
determined by its reproducing property, and holomorphic functions
on~\(G\) correspond to holomorphic functions on~\(\Omega\) that
are independent of the extra variable. Consequently, zeroes of
weighted Bergman kernel functions give rise to zeroes of ordinary
Bergman kernel functions of higher-dimensional domains.

As a concrete example, consider on a bounded domain~\(G\)
in~\(\C^n\) containing the origin the weight function
\(\exp(-t\|z\|)\), where \(t\) is a real parameter, and
\(\|z\|\)~denotes the Euclidean length
\(\sqrt{|z_1|^2+\dots+|z_n|^2}\).  Let \(K_t(z,w)\) denote the
Bergman kernel function of~\(G\) with respect to this weight. It
is a special case of a recent theorem of Miroslav Engli{\v{s}}
\cite{englis} that this weighted Bergman kernel function must
have zeroes near the origin when \(t\)~is sufficiently large.

Remarkably, the proof depends on the non-smoothness of the weight
function.  The idea is to show that
\(\lim_{t\to\infty}K_t(z,z)^{1/t}=e^{\|z\|}\). It follows that
the function \(K_t(z,w)\) cannot be zero-free for every
large~$t$, for if it were, then a sesqui-holomorphic branch of
\(K_t(z,w)^{1/t}\) could be defined near the origin.  When
\(t\to\infty\), there would be a limiting sesqui-holomorphic
function \(L(z,w)\) such that \(L(z,z)=e^{\|z\|}\); but this is
impossible because the function~\(e^{\|z\|}\) is not real
analytic at the origin.

The verification that
\(\lim_{t\to\infty}K_t(z,z)^{1/t}=e^{\|z\|}\) is carried out via
an upper estimate and a lower estimate.  Let \(\|f\|_t\) denote
the weighted norm \((\int_G |f(w)|^2 e^{-t\|w\|}\, dV_w)^{1/2}\).
If \(f\)~is a holomorphic function, and \(B_z\)~is a small ball
centered at~\(z\) with volume \(|B_z|\), then the mean-value
property of holomorphic functions and the Cauchy-Schwarz
inequality imply that \(|f(z)|\) is bounded above by \(\|f\|_{t}
\,|B_z|^{-1} (\int_{B_z} e^{t\|w\|}\, dV_w)^{1/2}\). Therefore
\(K_t(z,z)\), which is the square of the norm of the point
evaluation functional, is bounded above by \(|B_z|^{-1}
\sup_{w\in B_z} e^{t\|w\|}\), and so \(\limsup_{t\to \infty}
K_t(z,z)^{1/t} \le \sup_{w\in B_z} e^{\|w\|}\). Now let the
radius of the ball~\(B_z\) shrink to zero to conclude that
\(\limsup_{t\to \infty} K_t(z,z)^{1/t} \le e^{\|z\|}\).

For the lower bound, use the convexity of the Euclidean norm and
the representation of a supporting hyperplane as the zero set of
the real part of a linear holomorphic function. For each
point~\(z\) in the domain, there is a holomorphic function~\(g\)
such that \(\re g(w)\le\|w\|\) for all~\(w\), and \(\re
g(z)=\|z\|\).  Since \(K_t(z,z)\) is the square of the norm of
the point evaluation functional, it is no smaller than
\(|e^{tg(z)}|/\|e^{tg/2}\|_t^2\), which in turn is no smaller
than \(e^{t\|z\|}\) divided by the volume of the domain.
Consequently, \(\liminf_{t\to\infty} K_t(z,z)^{1/t}\ge
e^{\|z\|}\).

  \begin{problem}
    For concrete examples, determine how large \(t\) must be
    taken in Engli{\v{s}}'s theorem to guarantee that the
    weighted Bergman kernel function has zeroes.
  \end{problem}

\section{Weighted disk kernels and convex domains}
In the preceding section, I remarked that the Bergman kernel
function of a Hartogs domain in~\(\C^{n+1}\) is related to a
weighted Bergman kernel function on the base domain in~\(\C^n\).
For the same reason, a multi-dimensional domain that is fibered
over a one-dimensional base has a Bergman kernel function that is
related to a weighted Bergman kernel function on the base.  In
this section, I shall discuss an interesting example of this
general principle.

The domain in~\(\C^n\) defined by the inequality
\begin{equation}
\label{eq:ex}
|z_1| + |z_2|^{2/p_2} + \dots + |z_n|^{2/p_n} <1
\end{equation}
has a Bergman kernel function whose restriction to the
\(z_1\)-axis is proportional to the weighted Bergman kernel
function for the unit disk \(\{z\in\C: |z|<1\}\) with weight
\((1-|z|)^{p_2+\dots+p_n}\).  Here the~\(p_j\) can be arbitrary
positive real numbers, and the proportionality constant is the
volume of the \((n-1)\)-dimensional domain defined by the
inequality
\begin{equation*}
|z_2|^{2/p_2} + \dots + |z_n|^{2/p_n} <1.
\end{equation*}

Accordingly, it is useful to compute explicitly the weighted
Bergman kernel function~\(K_q\) for the unit disk with weight
\((1-|z|)^q\), where \(q>0\).  The square of the norm of the
monomial \(z^k\) with weight factor \( (1-|z|)^q\) is
\begin{equation*}
\int_0^{2\pi}\int_0^1 r^{2k+1} (1-r)^q \,dr\,d\theta = 
2\pi \Beta(2k+2,q+1),
\end{equation*}
where \(\Beta\) is the Beta function defined in terms of the
Gamma function by \(\Beta(a,b)= \Gamma(a)
\Gamma(b)/\Gamma(a+b)\).  Consequently, the weighted Bergman
kernel function \(K_q(z,w)\) equals
\((2\pi)^{-1}\sum_{k=0}^\infty (z\bar w)^k/\Beta(2k+2,q+1)\). A
closed form expression for this series is most conveniently
written in terms of the squares of the variables:
\begin{equation}
\label{eq:q}
K_q(z^2, w^2)= \frac{(q+1)}{4\pi z\bar w} \left[ \frac{1}{(1-z\bar
w)^{q+2}} - \frac{1}{(1+z\bar w)^{q+2}} \right].
\end{equation}
The powers of $(1\pm z\bar w)$ are to be understood as principal
branches.  The validity of this closed form expression can be
verified by the binomial series expansion.

The explicit expression~\eqref{eq:q} implies that the weighted
Bergman kernel function~\(K_q\) has zeroes in the interior of the
unit disk if and only if \(q>2\). Indeed, taking limits
in~\eqref{eq:q} shows that \(K_q\)~is not equal to~\(0\) when
either coordinate is equal to~\(0\), so it is only necessary to
decide if \((1-t)^{q+2}=(1+t)^{q+2}\) for some non-zero~\(t\) in
the unit disk. Since the mapping \(t\mapsto (1+t)/(1-t)\) takes
the unit disk bijectively to the right half-plane, with the
origin going to the point~\(1\), following this mapping by the
mapping \(u\mapsto u^{q+2}\) produces a composite mapping that
takes some non-zero point of the open unit disk to the
point~\(1\) if and only if \(q>2\).

Consequently, the Bergman kernel function of the
domain~\eqref{eq:ex} is guaranteed to have zeroes if
\(p_2+\dots+p_n >2\). Moreover, this domain is geometrically
convex if no~\(p_j\) exceeds~\(2\).  Therefore, one can exhibit
many concrete examples of convex domains whose Bergman kernel
functions have zeroes \cite{boasfustraube}. Here are some:
\begin{equation}
\label{eq:convex}
\begin{gathered}
\{\,z\in\C^3: |z_1|+|z_2|+|z_3|<1\,\}, \\
\{\,z\in\C^3: |z_1|+|z_2|+|z_3|^2<1\,\}, \\
\{\,z\in\C^4: |z_1|+|z_2|^2+|z_3|^2+|z_4|^4<1\,\}.
\end{gathered}
\end{equation}

Using a different method, Peter Pflug and E.~H. Youssfi found
some other interesting examples of convex domains whose Bergman
kernel functions have zeroes.  Even though the ``minimal ball''
in~\(\C^n\) defined by the inequality
\begin{equation}
\label{eq:minimal}
  |z_1|^2+\dots+|z_n|^2 + |z_1^2+\dots+z_n^2| <1
\end{equation}
lacks multi-circular symmetry, its Bergman kernel function is
known explicitly \cite{opy}, and in \cite{pflugyoussfi} the
authors analyzed the explicit formula to see that the Bergman
kernel function of this domain has zeroes when \(n\ge4\).

Although the convex domains defined by \eqref{eq:convex}
and~\eqref{eq:minimal} do not have smooth boundaries, they can be
approximated from inside by smoothly bounded, strongly convex
domains.  From the method of variation of domains in
section~\ref{variation-domains}, it follows that when \(n\ge3\),
there exist smoothly bounded, strongly convex domains in~\(\C^n\)
whose Bergman kernel functions have zeroes.  Pflug and Youssfi
even showed in \cite{pflugyoussfi} how to write down concrete
examples of bounded, smooth, algebraic, strongly convex domains
that approximate~\eqref{eq:minimal} from inside.

Using the same idea, Nguy\^en Vi\^et Anh \cite{anh} gave concrete
examples of bounded, smooth, algebraic, strongly convex,
Reinhardt domains in~\(\C^n\) whose Bergman kernel functions have
zeroes when \(n\ge3\). For example, when \(k\)~is a sufficiently
large positive integer, the inequality
\begin{equation}
\label{eq:alg}
(|z_1|^2+|z_2|^2+|z_3|^2)^{2k} + \sum_{\substack{ \pm \\ \text{(8
terms)}}} (\pm |z_1| \pm |z_2| \pm |z_3|)^{2k} <1
\end{equation}
defines such a domain in~\(\C^3\).  

To see that~\eqref{eq:alg} has the required properties, first
observe that the \(\ell_{2k}\) norm decreases to the
\(\ell_\infty\) norm as \(k\to\infty\), so these domains are
interior approximations to the domain \(\{\,z\in\C^3:
|z_1|+|z_2|+|z_3|<1\,\}\), which is one of the domains
\eqref{eq:convex} whose Bergman kernel functions have zeroes.
The odd powers of the \(|z_j|\) in the expansion
of~\eqref{eq:alg} cancel out by symmetry, so the defining
function is equivalent to a polynomial.  It would be obvious that
the defining function~\eqref{eq:alg} is convex if it had \(\re
z_j\) in place of~\(|z_j|\), for a convex function of a linear
function is convex; now observe that positive combinations of
even powers are increasing, and the composite of a convex
increasing function with the convex function~\(|z_j|\) is convex.

The two-dimensional domain defined by the inequality
\(|z_1|+|z_2|<1\) is a borderline case for the preceding
considerations. It turns out \cite{boasfustraube} that the
Bergman kernel function of this domain has no zeroes in the
interior of the domain, although it does have zeroes on the
boundary.

\begin{problem}
  Exhibit a bounded convex domain in~\(\C^2\) whose Bergman
  kernel function has zeroes in the interior of the domain.
\end{problem}

\section{Conclusion}
It is a difficult problem to determine whether the Bergman kernel
function of a specific domain has zeroes or not. If the kernel
function is presented as an infinite series, then locating the
zeroes may be of the same order of difficulty as proving the
Riemann hypothesis; and even if the series can be summed in
closed form, determining whether or not \(0\)~is in the range may
be hard.

In this article, I have emphasized examples in which the Bergman
kernel function does have zeroes. As the subject developed
historically, such examples were considered surprising. From our
current perspective, it would be more surprising to find some
simple geometric condition guaranteeing that the Bergman kernel
function is zero-free.

Students planning further investigation of the Bergman kernel
function might consult, in addition to the journal articles I
have cited, Bell's book \cite{bellbook} about the one-dimensional
theory, the book of Jarnicki and Pflug \cite{jarnickipflug}, and
Stefan Bergman's own book \cite{bergman}.  I offer the following
problem as an illustration of how much remains to be discovered
about the zeroes of the Bergman kernel function.

\begin{problem}
  Characterize the vectors \((p_1, p_2, \dots, p_n)\) of positive
  numbers for which the Bergman kernel function of the domain
  in~\(\C^n\) defined by the inequality
\begin{equation*}
|z_1|^{2/p_1}+|z_2|^{2/p_2} + \dots + |z_n|^{2/p_n} <1
\end{equation*}
is zero-free.
\end{problem}


\begin{thebibliography}{10}

\bibitem{anh} Nguy\^en~Vi\^et Anh, \emph{The {L}u {Q}i-{K}eng
conjecture fails for strongly convex algebraic complete
{R}einhardt domains in $\mathbf{C}\sp n\ (n\ge3)$},
Proc. Amer. Math. Soc., to appear.

\bibitem{belldavid} David Bell, \emph{Some properties of the
{B}ergman kernel function}, Compositio Math. \textbf{21} (1969),
329--330.

\bibitem{bell1981} Steven~R. Bell, \emph{Proper holomorphic
mappings and the {B}ergman projection}, Duke Math. J. \textbf{48}
(1981), no.~1, 167--175.

\bibitem{bell1982} \bysame, \emph{The {B}ergman kernel function
and proper holomorphic mappings},
Trans. Amer. Math. Soc. \textbf{270} (1982), no.~2, 685--691.

\bibitem{bellbook} \bysame, \emph{The {C}auchy transform,
potential theory, and conformal mapping}, CRC Press, 1992.

\bibitem{bergman} Stefan Bergman, \emph{The kernel function and
conformal mapping}, revised ed., American Mathematical Society,
1970.

\bibitem{boascounterexample} Harold~P. Boas, \emph{Counterexample
to the {L}u {Q}i-{K}eng conjecture}, Proc.
Amer. Math. Soc. \textbf{97} (1986), no.~2, 374--375.

\bibitem{boasgeneric} \bysame, \emph{The {L}u {Q}i-{K}eng
conjecture fails generically}, Proc. Amer.
Math. Soc. \textbf{124} (1996), no.~7, 2021--2027.

\bibitem{boasfustraube} Harold~P. Boas, Siqi Fu, and
Emil~J. Straube, \emph{The {B}ergman kernel function: explicit
formulas and zeroes}, Proc. Amer. Math. Soc. \textbf{127} (1999),
no.~3, 805--811.

\bibitem{chernji} Shiing-Shen Chern and Shanyu Ji, \emph{On the
{R}iemann mapping theorem}, Ann.  of Math. (2) \textbf{144}
(1996), no.~2, 421--439.

\bibitem{englis} Miroslav Engli{\v{s}}, \emph{Asymptotic
behaviour of reproducing kernels of weighted {B}ergman spaces},
Trans. Amer. Math. Soc. \textbf{349} (1997), no.~9, 3717--3735.

\bibitem{greenekrantz} R.~E. Greene and Steven~G. Krantz,
  \emph{Stability properties of the {B}ergman kernel and
    curvature properties of bounded domains}, Recent developments
  in several complex variables, Princeton Univ. Press, 1981,
  pp.~179--198.
  
\bibitem{jarnickipflug} Marek Jarnicki and Peter Pflug,
  \emph{Invariant distances and metrics in complex analysis},
  Walter de Gruyter \& Co., 1993.

\bibitem{kakurai}
Shigeki Kakurai, \emph{On the Lu Qi-Keng conjecture},
Rep. Fac. Engrg. Kanagawa Univ. no.~17 (1979), 3--4.

\bibitem{kanemaru}
Tadayoshi Kanemaru,  
\emph{A remark on the Lu Qi Keng conjecture},
Mem. Fac. Ed. Kumamoto Univ. Natur. Sci. no.~31 (1982), 1--3. 

\bibitem{Krantz} Steven~G. Krantz, \emph{Function theory of
several complex variables}, second ed., Wadsworth \& Brooks/Cole,
1992.

\bibitem{lu} Lu~Qi-keng, \emph{On {K}aehler manifolds with
constant curvature}, Chinese Math.--Acta \textbf{8} (1966),
283--298.

\bibitem{matsuura} Shozo Matsuura, \emph{On the {L}u {Q}i-{K}eng
conjecture and the {B}ergman representative domains}, Pacific
J. Math. \textbf{49} (1973), 407--416.

\bibitem{opy} K.~Oeljeklaus, P.~Pflug, and E.~H. Youssfi,
\emph{The {B}ergman kernel of the minimal ball and applications},
Ann. Inst. Fourier (Grenoble) \textbf{47} (1997), no.~3,
915--928.

\bibitem{pflugyoussfi} P.~Pflug and E.~H. Youssfi, \emph{The {L}u
{Q}i-{K}eng conjecture fails for strongly convex algebraic
domains}, Arch. Math. (Basel) \textbf{71} (1998), no.~3,
240--245.

\bibitem{ramadanov} I.~Ramadanov, \emph{Sur une propri\'et\'e de
la fonction de {B}ergman}, C. R.  Acad. Bulgare Sci. \textbf{20}
(1967), 759--762.

\bibitem{rosenthal} Paul Rosenthal, \emph{On the zeros of the
{B}ergman function in doubly-connected domains},
Proc. Amer. Math. Soc. \textbf{21} (1969), 33--35.

\bibitem{semmes} Stephen Semmes, \emph{A generalization of
{R}iemann mappings and geometric structures on a space of domains
in $\mathbf{C}\sp n$}, Mem. Amer. Math. Soc.  \textbf{98} (1992),
no.~472.

\bibitem{skwar} M.~Skwarczy{\'n}ski, \emph{The invariant distance
in the theory of pseudoconformal transformations and the {L}u
{Q}i-keng conjecture}, Proc.  Amer. Math. Soc. \textbf{22}
(1969), 305--310.

\bibitem{skwarczynski} Maciej Skwarczy{\'n}ski,
\emph{Biholomorphic invariants related to the {B}ergman
function}, Dissertationes Math. \textbf{173} (1980).

\bibitem{suitayamada} Nobuyuki Suita and Akira Yamada, \emph{On
the {L}u {Q}i-keng conjecture}, Proc.
Amer. Math. Soc. \textbf{59} (1976), no.~2, 222--224.

\end{thebibliography}
\end{document}